\documentclass[a4paper,11pt]{amsart}
\usepackage{graphicx}
 \usepackage{mathptmx}
\usepackage{amsmath}
\usepackage{amssymb}
\usepackage{enumitem}
\usepackage{xcolor}

\newmuskip\pFqmuskip

\newcommand*\pFq[6][8]{%
  \begingroup % only local assignments
  \pFqmuskip=#1mu\relax
  \mathcode`=\string"8000
  \begingroup\lccode`\~=`\,
  \lowercase{\endgroup\let~}\pFqcomma
  F^{#2}_{#3}{\left(\genfrac..{0pt}{}{#4}{#5}\bigg|#6\right)}%
  \endgroup
}
\newcommand{\pFqcomma}{\mskip\pFqmuskip}

\newtheorem{theorem}{Theorem}
\newtheorem{lemma}[theorem]{Lemma}
\newtheorem{corollary}[theorem]{Corollary}

\begin{document}

\title[Degenerate polyexponential functions and degenerate Bell polynomials]{Degenerate polyexponential functions and degenerate Bell polynomials}

\author{Taekyun  Kim}
\address{Department of Mathematics, Kwangwoon University, Seoul 139-701, Republic of Korea}
\email{tkkim@kw.ac.kr}

\author{Dae San Kim}
\address{Department of Mathematics, Sogang University, Seoul, Republic of Korea}
\email{dskim@sogang.ac.kr}

\subjclass[2010]{11B73; 11B83; 05A19}
\keywords{degenerate polyexponential function; new type degenerate Bell polynomial; degenerate Lerch zeta function; degenerate Hurwitz zeta function; degenerate Riemann zeta function}

\maketitle

\begin{abstract}
In recent years, studying degenerate versions of some special polynomials, which was initiated by Carlitz in an investigation of the degenerate Bernoulli and Euler polynomials, regained lively interest of many mathematicians. 
In this paper, as a degenerate version of polyexponential functions introduced by Hardy, we study degenerate polyexponential functions and derive various properties of them. Also, we introduce new type degenerate Bell polynomials, which are different from the previously studied partially degenerate Bell polynomials and arise naturally in the recent study of degenerate zero-truncated Poisson random variables, and deduce some of their properties. Furthermore, we derive some identities connecting the polyexponential functions and the new type degenerate Bell polynomials.
\end{abstract}

\section{Introduction}
It is well known that the Stirling numbers of the second kind are defined by
\begin{equation}
    x^{n}=\sum_{l=0}^{n}S_{2}(n,l)(x)_{l},\quad (n\ge 0),\quad (\mathrm{see}\ [4,6,8,12,22]),\label{1}
\end{equation}
where $(x)_{0}=1$, $(x)_{n}=x(x-1)\cdots(x-n+1)$, $(n\ge 1)$.\\
From \eqref{1}, we note that
\begin{equation}
    \frac{1}{k!}(e^{t}-1)^{k}=\sum_{n=k}^{\infty}S_{2}(n,k)\frac{t^{n}}{n!},\quad(k\ge 0), (\mathrm{see}\ [4,5,7,16]).\label{2}
\end{equation}
The Bell polynomials are given by the generating function 
\begin{equation}
    e^{x(e^{t}-1)}=\sum_{l=0}^{\infty}\mathrm{Bel}_{l}(x)\frac{t^{l}}{k!},\quad (\mathrm{see}\ [4,13]).\label{3}
\end{equation}
When $x=1$, $\mathrm{Bel}_{n}=\mathrm{Bel}_{n}(1)$ is called the $n$-th Bell number. \\
From \eqref{3}, we note that
\begin{equation}
    \mathrm{Bel}_{n}(x)=\sum_{k=0}^{n}x^{k}S_{2}(n,k),\quad (n\ge 0),\quad (\mathrm{see}\ [4,13]).\label{4}
\end{equation}
For $k\in\mathbb{N}\cup\{0\}$, the polylogarithms are defined as
\begin{equation}
    \mathrm{Li}_{k}(x)=\sum_{n=1}^{\infty}\frac{x^{n}}{n^{k}},\quad (\mathrm{see}\ [4,12]).\label{5}
\end{equation}
Note that $\mathrm{Li}_{1}(x)=-\log(1-x)$. \\
The polyexponential function was first considered by Hardy and is given by
\begin{equation}
    e(x,\ a|s)=\sum_{n=0}^{\infty}\frac{x^{n}}{(n+a)^{s}n!},\quad(\mathrm{Re}(a)>0),\quad (\mathrm{see}\ [9,10,11]).\label{6}
\end{equation}
From \eqref{6}, we note that
\begin{equation}
    e(x,\ 1|1)=\sum_{n=0}^{\infty}\frac{x^{n}}{(n+1)!}=\frac{1}{x}\bigg(\sum_{n=0}^{\infty}\frac{x^{n}}{n!}-1\bigg)=\frac{1}{x}(e^{x}-1),\quad (\mathrm{see}\ [2,9,10,11]).\label{7}
\end{equation}
~~\\
For any nonzero $\lambda\in\mathbb{R}$, the degenerate exponential function is defined by
\begin{equation}
    e_{\lambda}^{x}(t)=(1+\lambda t)^{\frac{x}{\lambda}},\quad e_{\lambda}(t)= e_{\lambda}^{1}(t),\quad(\mathrm{see}\ [14,15,19]).\label{8}
\end{equation}
Note that $\displaystyle \lim_{\lambda\rightarrow 0}e^{x}_{\lambda}(t)=e^{xt}$. \\
In [3], Carlitz considered the degenerate Bernoulli polynomials which are given by 
\begin{equation}
    \frac{t}{e_{\lambda}(t)-1}e_{\lambda}^{x}(t)=\sum_{n=0}^{\infty}\beta_{n,\lambda}(x)\frac{t^{n}}{n!}.\label{9}
\end{equation}
When $x=0$, $\beta_{n,\lambda}=\beta_{n,\lambda}(0)$ is called the $n$-th degenerate Bernoulli number.
Note that $\displaystyle\lim_{\lambda\rightarrow 0}\beta_{n,\lambda}(x)=B_{n}(x)\displaystyle$, where $B_{n}(x)$ are the ordinary Bernoulli polynomials given by
\begin{displaymath}
    \frac{t}{e^{t}-1}e^{xt}=\sum_{n=0}^{\infty}B_{n}(x)\frac{t^{n}}{n!},\quad(\mathrm{see}\ [12-22]).
\end{displaymath}
Recently, Kim introduced the degenerate Stirling numbers of the second kind which are defined by
\begin{equation}
    (x)_{n,\lambda}=\sum_{k=0}^{n}S_{2,\lambda}(n,k)(x)_{k},\quad (n\ge 0),\quad(\mathrm{see}\ [16]),\label{10}
\end{equation}
where $(x)_{0,\lambda}=1$, $(x)_{b,\lambda}=x(x-\lambda)(x-2\lambda)\cdots(x-(n-1)\lambda)$, $(n\ge 1)$. \\
In [13], the partially degenerate Bell polynomial are defined as
\begin{equation}
    e^{x(e_{\lambda}(t)-1)}=\sum_{n=0}^{\infty}\mathrm{bel}_{n,\lambda}(x)\frac{t^{n}}{n!}. \label{11}
\end{equation}
From \eqref{10} and \eqref{11}, we note that
\begin{equation}
    \mathrm{bel}_{n,\lambda}(x)=\sum_{k=0}^{n}x^{k}S_{2,\lambda}(n,k),\quad(\mathrm{see}\ [13,20]).\label{12}
\end{equation}
It is easy to show that $\displaystyle\lim_{\lambda\rightarrow 0}\mathrm{bel}_{n,\lambda}(x)=\mathrm{Bel}_{n}(x),\ (n\ge 0)\displaystyle$. \\
By \eqref{11}, we easily get
\begin{equation}
    \mathrm{bel}_{n,\lambda}(x)=e^{-x}\sum_{k=0}^{\infty}\frac{(k)_{n,\lambda}}{k!}x^{k},\quad (n\ge 0).\label{13}
\end{equation}

In this paper, as a degenerate version of polyexponential functions introduced by Hardy, we study degenerate polyexponential functions and derive various properties of them. Also, we introduce new type degenerate Bell polynomials, which are different from the previously studied partially degenerate Bell polynomials and arise naturally in the recent study of degenerate zero-truncated Poisson random variables, and deduce some of their properties. Furthermore, we derive some identities connecting the polyexponential functions and the new type degenerate Bell polynomials. \\
\indent This paper is organized as follows. In Section 1, we recall some necessary ingredients, namely polyexponential functions, degenerate Bernoulli polynomials, degenerate Stirling numbers of the second kind and partially degenerate Bell polynomials. In Section 2, we introduce the degenerate polyexponential functions and derive several properties of them. Especially, we find integral representations for the degenerate polyexponential functions. Then we introduce the new type degenerate Bell polynomials, deduce some of their properties and get their connections with degenerate polyexponential functions. Finally, in Section 3, we recall the degenerate Lerch zeta functions which are closely related to the degenerate polyexponential functions. Also, we note the degenerate Hurwitz zeta functions and the degenerate Riemann zeta functions, as special cases of the degenerate Lerch zeta functions.

\section{Degenerate polyexponential functions and degenerate Bell polynomials}
In view of \eqref{8}, we consider the degenerate polyexponential functions which are given by
\begin{equation}
    e_{\lambda}(x,\ \delta|k)=\sum_{n=0}^{\infty}\frac{(1)_{n,\lambda}}{n!(n+\delta)^{k}}x^{n},\label{14}
\end{equation}
where $k\in\mathbb{N}\cup\{0\}$, and $\delta\in\mathbb{C}$ with $\mathrm{Re}(\delta)>0$. \\
By \eqref{14}, we easily get
\begin{equation}
    e_{\lambda}(x,\ \delta|0)=e_{\lambda}(x),\quad e_{\lambda}(x,\ 1|1)=\frac{1}{x}\frac{1}{1+\lambda}(e_{\lambda}(x)-1)+\frac{\lambda}{1+\lambda}e_{\lambda}(x).\label{15}
\end{equation}
Note that
\begin{displaymath}
    \lim_{\lambda\rightarrow 0}e_{\lambda}(x,\ 1|1)=\frac{1}{x}(e^{x}-1)=e(x,\ 1|1).
\end{displaymath}
It is known that the degenerate gamma function is defined by
\begin{equation}
    \Gamma_{\lambda}(s)=\int_{0}^{\infty}t^{s-1}e_{\lambda}^{-1}(t)dt,\quad(\mathrm{Re}(s)>0),\quad (\mathrm{see}\ [15,19]).\label{16}
\end{equation}
From \eqref{16}, we consider the incomplete degenerate gamma function given by
\begin{equation}
    d_{\lambda}(\delta,x)=\int_{0}^{x}t^{\delta-1}e_{\lambda}(-t)dt=x^{\delta}e_{\lambda}(-x,\ \delta|1),\label{17}
\end{equation}
where $\mathrm{Re}(\delta)>0$ and $x\ge 0$.
Note that $\displaystyle\lim_{x\rightarrow \infty}d_{\lambda}(\delta,\ x)=\Gamma_{-\lambda}(\delta)\displaystyle$.
\begin{theorem}
    For $\delta\in\mathbb{C}$ with $\mathrm{Re}(\delta)>0$, and $x\ge 0$, we have
    \begin{displaymath}
        x^{\delta}e_{\lambda}(-x,\ \delta|1)=d_{\lambda}(\delta, x).
    \end{displaymath}
\end{theorem}
From \eqref{17}, we note that
\begin{equation}
    \frac{d}{dx}[x^{\delta}e_{\lambda}(-x,\ \delta|1)]=x^{\delta-1}e_{\lambda}(-x). \label{18}
\end{equation}
In particular,
\begin{displaymath}
    e_{\lambda}(x,\ 1|2)=\sum_{n=0}^{\infty}\frac{(1)_{n,\lambda}}{n!(n+1)^{2}}x^{n},
\end{displaymath}
and
\begin{align}
    \frac{1}{x}\int_{0}^{x}\frac{1}{t}(e_{\lambda}(t)-1)dt&=\frac{1}{x}\int_{0}^{x}\frac{1}{t}\bigg(\sum_{n=1}^{\infty}\frac{(1)_{n,\lambda}}{n!}t^{n}\bigg)dt=\frac{1}{x}\sum_{n=0}^{\infty}\frac{(1)_{n+1,\lambda}}{(n+1)!}\int_{0}^{x}t^{n}dt\label{19}\\
    &=\frac{1}{x}\sum_{n=0}^{\infty}\frac{(1)_{n,\lambda}}{(n+1)!(n+1)}(1-n\lambda)x^{n+1}\nonumber\\
    &=\sum_{n=0}^{\infty}\frac{(1)_{n,\lambda}}{n!(n+1)^{2}}x^{n}-\sum_{n=0}^{\infty}\frac{(1)_{n,\lambda}(n+1-1)\lambda}{n!(n+1)^{2}}x^{n}\nonumber \\
    &=e_{\lambda}(x,\ 1|2)-\lambda\big(e_{\lambda}(x,\ 1|1)-e_{\lambda}(x,\ 1|2)\big)\nonumber \\
    &=(1+\lambda)e_{\lambda}(x,\ 1|2)-\lambda e_{\lambda}(x,\ 1|1). \nonumber
\end{align}
By \eqref{15} and \eqref{19}, we get
\begin{align}
    e_{\lambda}(x,\ 1|2)&=\frac{1}{1+\lambda}\cdot\frac{1}{x}\int_{0}^{x}\frac{1}{t}\big(e_{\lambda}(t)-1\big)dt+\frac{\lambda}{1+\lambda}e_{\lambda}(x,\ 1|1) \label{20}\\
    &=\frac{1}{1+\lambda}\cdot\frac{1}{x}\int_{0}^{x}\frac{1}{t}(e_{\lambda}(t)-1)dt+\frac{1}{x}\frac{\lambda}{(1+\lambda)^{2}}(e_{\lambda}(x)-1)+\bigg(\frac{\lambda}{1+\lambda}\bigg)^{2}e_{\lambda}(x). \nonumber
\end{align}
\begin{lemma}
    For all $x$, we have
    \begin{displaymath}
        e_{\lambda}(x,\ 1|2)=\frac{1}{1+\lambda}\cdot\frac{1}{x}\int_{0}^{x}\frac{1}{t}(e_{\lambda}(t)-1)dt+\frac{1}{x}\frac{\lambda}{(1+\lambda)^{2}}(e_{\lambda}(x)-1)+\bigg(\frac{\lambda}{1+\lambda}\bigg)^{2}e_{\lambda}(x).
    \end{displaymath}
\end{lemma}
Let $\mathrm{Ein}_{\lambda}(z)$ be the special  degenerate entire function which is defined by
\begin{align}
    \mathrm{Ein}_{\lambda}(x)&=\int_{0}^{x}\frac{1}{t}\big(1-e_{\lambda}(-t)\big)dt=\int_{0}^{x}\frac{1}{t}\bigg(1-\sum_{k=0}^{\infty}\frac{(1)_{k,\lambda}}{k!}(-1)^{k}t^{k}\bigg)dt \label{21} \\
    &=\int_{0}^{x}\frac{1}{t}\sum_{k=1}^{\infty}(-1)^{k-1}\frac{(1)_{k,\lambda}}{k!}t^{k}dt=\sum_{k=1}^{\infty}\frac{(1)_{k,\lambda}}{k!}(-1)^{k-1}\int_{0}^{x}t^{k-1}dt \nonumber \\
    &=\sum_{k=1}^{\infty}\frac{(-1)^{k-1}}{k!\cdot k}(1)_{k,\lambda}x^{k}. \nonumber
\end{align}
On the other hand,
\begin{align}
    e_{\lambda}(x,\ 1|2)&=\frac{1}{1+\lambda}\frac{1}{x}\int_{0}^{x}\frac{1}{t}\big(e_{\lambda}(t)-1\big)dt+\frac{1}{x}\frac{\lambda}{(1+\lambda)^{2}}(e_{\lambda}(x)-1)+\bigg(\frac{\lambda}{1+\lambda}\bigg)^{2}e_{\lambda}(x)\label{22}\\
    &=-\frac{1}{1+\lambda}\frac{1}{x}\int_{0}^{-x}\frac{1}{t}\big(1-e_{\lambda}(-t)\big)dt+\frac{1}{x}\frac{\lambda}{(1+\lambda)^{2}}(e_{\lambda}(x)-1)+\bigg(\frac{\lambda}{1+\lambda}\bigg)^{2}e_{\lambda}(x) \nonumber\\
    &=-\frac{1}{1+\lambda}\frac{1}{x}\mathrm{Ein}_{\lambda}(-x)+\frac{1}{x}\frac{\lambda}{(1+\lambda)^{2}}(e_{\lambda}(x)-1)+\bigg(\frac{\lambda}{1+\lambda}\bigg)^{2}e_{\lambda}(x). \nonumber
\end{align}
Thus, we have
\begin{equation}
    x(1+\lambda)e_{\lambda}(-x,\ 1|2)=\mathrm{Ein}_{\lambda}(x)-\frac{\lambda}{1+\lambda}\big(e_{\lambda}(-x)-1\big)+x\frac{\lambda^{2}}{1+\lambda}e_{\lambda}(-x).\label{23}
\end{equation}
From \eqref{23}, we note that
\begin{displaymath}
    \mathrm{Ein}_{\lambda}(x)=x(1+\lambda)e_{\lambda}(-x,\ 1|2)+\frac{\lambda}{1+\lambda}\big(e_{\lambda}(-x)-1\big)-x\frac{\lambda^{2}}{1+\lambda}e_{\lambda}(-x).
\end{displaymath}
It is easy to show that
\begin{displaymath}
    \bigg(\frac{d}{dx}x\bigg)e_{\lambda}(x,\ 1|p+1)=\frac{d}{dx}x\sum_{n=0}^{\infty}\frac{(1)_{n,\lambda}x^{n}}{n!(n+1)^{p+1}}=e_{\lambda}(x,\ 1|p).
\end{displaymath}
Thus, we have
\begin{displaymath}
    e_{\lambda}(x,\ 1|p+1)=\frac{1}{x}\int_{0}^{x}e_{\lambda}(t,\ 1|p)dt.
\end{displaymath}
In general,
\begin{displaymath}
    \bigg(\frac{d}{dx}x\bigg)\big[x^{\delta-1}e_{\lambda}(x,\ \delta|p+1)\big]=\frac{d}{dx}\sum_{n=0}^{\infty}\frac{(1)_{n,\lambda}x^{n+\delta}}{n!(n+\delta)^{p+1}}=x^{\delta-1}e_{\lambda}(x,\ \delta|p).
\end{displaymath}
Thus, we have
\begin{displaymath}
    e_{\lambda}(x,\ \delta|p+1)=\frac{1}{x^{\delta}}\int_{0}^{x}t^{\delta-1}e_{\lambda}(x,\ \delta|p)dt.
\end{displaymath}
For $p\in\mathbb{N}$, we have
\begin{align}
    \bigg(\frac{d}{dx}x\bigg)^{p}\big[x^{\delta-1}e_{\lambda}(x,\ \delta|p)\big]&=\underbrace{\bigg(\frac{d}{dx}x\bigg)\times \bigg(\frac{d}{dx}x\bigg)\times\cdots\times\bigg(\frac{d}{dx}x\bigg)}_{p-\mathrm{times}}\big[x^{\delta-1}e_{\lambda}(x,\ \delta|p)\big]\label{24} \\
    &=\bigg(\frac{d}{dx}x\bigg)^{p-1}\big[x^{\delta-1}e_{\lambda}(x,\ \delta|p-1)\big]\nonumber \\
    &=\cdots\nonumber\\
    &=x^{\delta-1}\sum_{n=0}^{\infty}\frac{(1)_{n,\lambda}}{n!}x^{n}=x^{\delta-1}e_{\lambda}(x).\nonumber
\end{align}
Let $|z|<|\delta|$. Then we note that
\begin{align}
    \sum_{p=0}^{\infty}e_{\lambda}(x,\ \delta|p)z^{p}&=\sum_{p=0}^{\infty}\sum_{n=0}^{\infty}\frac{(1)_{n,\lambda}x^{n}}{n!(n+\delta)^{p}}z^{p}=\sum_{n=0}^{\infty}\frac{(1)_{n,\lambda}}{n!}x^{n}\sum_{p=0}^{\infty}\bigg(\frac{z}{n+\delta}\bigg)^{p}\nonumber \\
    &=\sum_{n=0}^{\infty}\frac{(1)_{n,\lambda}}{n!}x^{n}\frac{n+\delta}{n+\delta-z}=\sum_{n=0}^{\infty}\frac{(1)_{n,\lambda}}{n!}x^{n}\bigg(1+\frac{z}{n+\delta-z}\bigg)\label{25}\\
    &=\sum_{n=0}^{\infty}\frac{(1)_{n,\lambda}}{n!}x^{n}+z\sum_{n=0}^{\infty}\frac{(1)_{n,\lambda}}{n!(n+\delta-z)}x^{n}\nonumber\\
    &=e_{\lambda}(x)+z\,e_{\lambda}(x,\ \delta-z|1). \nonumber
\end{align}
Therefore, we obtain the following theorem.
\begin{theorem}
    For $p\in\mathbb{N}$ and $|z|<|\delta|$, we have
    \begin{displaymath}
        \bigg(\frac{d}{dx}x\bigg)^{p}\big[x^{\delta-1}e_{\lambda}(x,\ \delta|p)\big]=x^{\delta-1}e_{\lambda}(x),
    \end{displaymath}
    and
    \begin{displaymath}
        \sum_{k=0}^{\infty}e_{\lambda}(x,\ \delta|k)z^{k}=e_{\lambda}(x)+z\,e_{\lambda}(x,\ \delta-z|1).
    \end{displaymath}
\end{theorem}
Now, we observe that
\begin{align}
    \bigg(\frac{d}{dz}\bigg)^{m}\big(e_{\lambda}(x,\ \delta-z|1)\big)&=\bigg(\frac{d}{dz}\bigg)^{m}\sum_{n=0}^{\infty}\frac{(1)_{n,\lambda}x^{n}}{n!(n+\delta-z)}=\sum_{n=0}^{\infty}\frac{(1)_{n,\lambda}x^{n}m!}{n!(n+\delta-z)^{m+1}}\label{26} \\
    &=m!e_{\lambda}(x,\ \delta-z|m+1). \nonumber
\end{align}
By using the Taylor expansion of $e_{\lambda}(x,\ \delta-z|1)$ with respect to $z$, we get
\begin{equation}
    e_{\lambda}(x,\ \delta-z|1)=\sum_{m=0}^{\infty}\bigg(\frac{d}{dz}\bigg)^{m}\big(e_{\lambda}(x,\ \delta-z|1\big)\bigg|_{z=0}\cdot\frac{z^{m}}{m!}=\sum_{m=0}^{\infty}e_{\lambda}(x,\ \delta|m+1)z^{m}. \label{27}
\end{equation}
Moreover,
\begin{equation}\begin{split}
  \bigg(\frac{d}{dz}\bigg)^{m}\big(e_{\lambda}(x,\ \delta-z|s)\big)&=\bigg(\frac{d}{dz}\bigg)^{m}\sum_{n=0}^{\infty}\frac{(1)_{n,\lambda}x^{n}}{n!(n+\delta-z)^{s}}. \label{28} \\
    &=\sum_{n=0}^{\infty}\frac{(1)_{n,\lambda}x^{n}}{n!}\cdot\frac{s(s+1)\cdots(s+m-1)}{(n+\delta-z)^{m+s}}\\ 
&=\frac{\Gamma(s+m)}{\Gamma(s)}e_{\lambda}(x,\ \delta-z|m+s).
\end{split}\end{equation}

By using Taylor expansion of $e_{\lambda}(x,\ \delta-z|s)$ with respect to $z$, we get
\begin{equation}
    e_{\lambda}(x,\ \delta-z|s)=\sum_{m=0}^{\infty}\frac{\big(\frac{d}{dz}\big)^{m}\big(e_{\lambda}(x,\ \delta-z|s)\big)\big|_{z=0}}{m!}z^{m}=\sum_{m=0}^{\infty}\frac{\Gamma(s+m)}{\Gamma(s)m!}e_{\lambda}(x,\ \delta|m+s)z^{m}.\label{29}
\end{equation}
Therefore, from \eqref{29}, we obtain the following theorem.
\begin{theorem}
    For $\delta,s\in\mathbb{C}$ with $\mathrm{Re}(\delta)>0,\ \mathrm{Re}(s)>0$, we have
    \begin{displaymath}
        e_{\lambda}(x,\ \delta-z|s)=\sum_{m=0}^{\infty}\frac{\Gamma(s+m)}{\Gamma(s)m!}e_{\lambda}(x,\ \delta|m+s)z^{m}.
    \end{displaymath}
    In addition, for $s=k\in\mathbb{N}$, we have
    \begin{displaymath}
        e_{\lambda}(x,\ \delta-z|k)=\sum_{m=0}^{\infty}\binom{k+m-1}{k-1}e_{\lambda}(x,\ \delta|m+k)z^{m}.
    \end{displaymath}
\end{theorem}
As is traditional, let $\int_{\infty}^{(0+)}$ (respectively, $\int_{-\infty}^{(0+)})$ denote the integration of path that starts at infinity (respectively, at negative infinity), encircles the origin counter-clockwise direction and returns to the starting point.
Then we have the following well-known integral representation of $-2i\sin\pi s\,\Gamma(s)$ as an entire function of $s$ (see \cite{21}):
\begin{equation}
\int_{\infty}^{(0+)}(-t)^{s-1}e^{-t}dt=-2i\sin\pi s\int_{0}^{\infty}t^{s-1}e^{-t}=-2i\sin\pi s\,\Gamma(s)=-2\pi i\frac{1}{\Gamma(1-s)}.\label{30}
\end{equation}
By \eqref{30}, we get
\begin{equation}
    \frac{1}{\Gamma(1-s)}=-\frac{1}{2\pi i}\int_{\infty}^{(0+)}(-t)^{s-1}e^{-t}dt=\frac{1}{2\pi i}\int_{-\infty}^{(0+)} t^{s-1}e^{t}dt.\label{31}
\end{equation}
From \eqref{31}, we easily note that
\begin{equation}\begin{split}\label{32}
\frac{\Gamma(1-s)}{2\pi i}\int_{-\infty}^{(0+)}t^{s-1}e^{\delta t}e^{nt}dt&=\frac{\Gamma(1-s)}{2\pi i}\frac{1}{(n+\delta)^{s}}\int_{-\infty}^{(0+)}e^{t}t^{s-1}dt \\
&=\frac{\Gamma(1-s)}{2\pi i}\frac{1}{(n+\delta)^{s}}\frac{2\pi i}{\Gamma(1-s)}=\frac{1}{(n+\delta)^{s}}.\nonumber
\end{split}\end{equation}

Thus, by \eqref{14} and \eqref{32}, we get
\begin{align}
    e_{\lambda}(x,\ \delta|s)&=\sum_{n=0}^{\infty}\frac{(1)_{n,\lambda}}{n!}\frac{x^{n}}{(n+\delta)^{s}}=\sum_{n=0}^{\infty}\frac{(1)_{n,\lambda}}{n!}x^n\frac{\Gamma(1-s)}{2\pi i}\int_{-\infty}^{(0+)}t^{s-1}e^{\delta t}e^{nt}dt \label{33}\\
    &=\frac{\Gamma(1-s)}{2\pi i}\int_{-\infty}^{(0+)}e^{\delta t}t^{s-1}\sum_{n=0}^{\infty}\frac{(1)_{n,\lambda}}{n!}e^{nt}x^{n}dt\nonumber \\
    &=\frac{\Gamma(1-s)}{2\pi i}\int_{-\infty}^{(0+)}e^{\delta t}t^{s-1}e_{\lambda}(xe^{t})dt.\nonumber
\end{align}
For $m\in\mathbb{N}$, we have
\begin{align}
    e_{\lambda}(x,\ \delta|m )&=\lim_{s\rightarrow m}e_{\lambda}(x,\ \delta|s)=\lim_{s\rightarrow m}\frac{1}{2\pi i}\frac{\pi}{\Gamma(s)\sin\pi s}\int_{-\infty}^{(0+)}e^{\delta t}t^{s-1}e_{\lambda}(xe^{t})dt \label{34}\\
    &=\frac{1}{2\pi i}\frac{(-1)^{m}}{\Gamma(m)}\int_{-\infty}^{(0+)}e^{\delta t}t^{m-1}e_{\lambda}(xe^{t})\mathrm{Log}\ tdt \nonumber\\
    &=\frac{1}{2\pi i}\frac{(-1)^{m}}{(m-1)!}\int_{-\infty}^{(0+)}e^{\delta t}t^{m-1}e_{\lambda}(xe^{t})\mathrm{Log}\ tdt. \nonumber
\end{align}
Thefore, by \eqref{33} and \eqref{34}, we obtain the following theorem.
\begin{theorem}
    For $s,\delta\in\mathbb{C}$ with $\mathrm{Re}(\delta)>0$, $s\in\mathbb{C}\backslash\mathbb{N}$, and $x\in\mathbb{C}$, we have the integral representation 
    \begin{displaymath}
        e_{\lambda}(x,\ \delta|s)=\frac{\Gamma(1-s)}{2\pi i}\int_{-\infty}^{(0+)}t^{s-1}e^{\delta t}e_{\lambda}e(xe^{t})\,dt.
    \end{displaymath}
In addition, for $s=m\in\mathbb{N}$, we have the integral representation
\begin{displaymath}
    e_{\lambda}(x,\ \delta|m)=\frac{(-1)^{m}}{2\pi i(m-1)!}\int_{-\infty}^{(0+)}t^{m-1}e^{\delta t}e_{\lambda}(xe^{t})\mathrm{Log}\ t \,dt.
\end{displaymath}
\end{theorem}
Now, we consider new type degenerate Bell polynomials which are given by
\begin{equation}
    e_{\lambda}(xe^{t})\cdot e_{\lambda}^{-1}(x)=\sum_{n=0}^{\infty}\mathrm{Bel}_{n,\lambda}(x)\frac{t^{n}}{n!}. \label{35}
\end{equation}
From \eqref{3} and \eqref{35}, we note that $\displaystyle\lim_{\lambda\rightarrow 0}\mathrm{Bel}_{n,\lambda}(x)=\mathrm{Bel}_{n}(x),\ (n\ge 0)\displaystyle$. \\
Now, we observe that
\begin{equation}
  \bigg(x\frac{d}{dx}\bigg)^{n}e_{\lambda}(x)=\bigg(x\frac{d}{dx}\bigg)^{n}\sum_{k=0}^{\infty}\frac{(1)_{k,\lambda}}{k!}x^{k}=\sum_{k=0}^{\infty}\frac{(1)_{k,\lambda}}{k!}k^{n}x^{k}. \label{36}
\end{equation}
On the other hand,
\begin{equation}
\sum_{n=0}^{\infty}\mathrm{Bel}_{n,\lambda}(x)\frac{t^{n}}{n!}=e_{\lambda}^{-1}(x)e_{\lambda}(xe^{t})=\sum_{n=0}^{\infty}\bigg(e_{\lambda}^{-1}(x)\sum_{k=0}^{\infty}\frac{(1)_{k,\lambda}}{k!}k^{n}x^k\bigg)\frac{t^{n}}{n!}. \label{37}
\end{equation}
When $x=1$, $\mathrm{Bel}_{n,\lambda}=\mathrm{Bel}_{n,\lambda}(1)$ are called the new type degenerate Bell numbers. \\
Therefore, by \eqref{36} and \eqref{37}, we obtain the following theorem.
\begin{theorem}[Dobinski's formula]
    For $n\ge 0$, we have
    \begin{displaymath}
        \mathrm{Bel}_{n,\lambda}(x)=e_{\lambda}^{-1}(x)\bigg(x\frac{d}{dx}\bigg)^{n}e_{\lambda}(x)=e_{\lambda}^{-1}(x)\sum_{k=0}^{\infty}\frac{(1)_{k,\lambda}}{k!}k^{n}x^{k}.
    \end{displaymath}
\end{theorem}
From \eqref{13}, we note that
\begin{align}
    e_{\lambda}(x,\ \delta|s)&=\sum_{k=0}^{\infty}\frac{(1)_{k,\lambda}x^{k}}{k!(k+\delta)^{s}}=\sum_{k=0}^{\infty}\frac{(1)_{k,\lambda}x^{k}}{k!}\frac{1}{\delta^{s}}\bigg(1+\frac{k}{\delta}\bigg)^{-s} \label{38} \\
    &=\sum_{k=0}^{\infty}\frac{(1)_{k,\lambda}}{k!}x^{k}\delta^{-s}\sum_{n=0}^{\infty}\binom{-s}{n}\bigg(\frac{k}{\delta}\bigg)^{n}\nonumber\\
    &=e_{\lambda}(x)\sum_{n=0}^{\infty}\binom{-s}{n}\frac{1}{\delta^{s+n}}e_{\lambda}^{-1}(x)\sum_{k=0}^{\infty}\frac{k^{n}}{k!}(1)_{k,\lambda}x^{k} \nonumber \\
    &=e_{\lambda}(x)\sum_{n=0}^{\infty}\binom{-s}{n}\frac{1}{\delta^{s+n}}\mathrm{Bel}_{n,\lambda}(x). \nonumber
\end{align}
Therefore, by \eqref{38}, we obtain the following corollary.
\begin{corollary}
    For $n\ge 0$ and $\mathrm{Re}(s)>0$, $\mathrm{Re}(\delta)>0$, we have
    \begin{displaymath}
        e_{\lambda}(x,\ \delta|s)=e_{\lambda}(x)\sum_{n=0}^{\infty}\binom{-s}{n}\frac{1}{\delta^{s+n}}\mathrm{Bel}_{n,\lambda}(x).
    \end{displaymath}
\end{corollary}
From \eqref{35}, we note that
\begin{equation}
    \frac{d}{dt}\big(e_{\lambda}(xe^{t})e_{\lambda}^{-1}(x)\big)=xe^{t}e_{\lambda}(xe^{t})e_{\lambda}^{-1}(x)\cdot (1+\lambda xe^{t})^{-1}. \label{39}
\end{equation}
Thus we have
\begin{equation}
    (1+\lambda xe^{t})\frac{d}{dt}\big(e_{\lambda}(xe^{t})e_{\lambda}^{-1}(x)\big)=xe^{t}e_{\lambda}(xe^{t})e_{\lambda}^{-1}(x)=xe^{t}\sum_{l=0}^{\infty}\mathrm{Bel}_{l,\lambda}(x)\frac{t^{l}}{l!} \label{40}
\end{equation}
\begin{displaymath}
    =x\sum_{m=0}^{\infty}\frac{t^{m}}{m!}\sum_{l=0}^{\infty}\mathrm{Bel}_{l,\lambda}(x)\frac{t^{l}}{l!}=x\sum_{n=0}^{\infty}\bigg(\sum_{l=0}^{n}\binom{n}{l}\mathrm{Bel}_{l,\lambda}(x)\bigg)\frac{t^{n}}{n!}.\qquad
\end{displaymath}
On the other hand,
\begin{align}
    (1+\lambda xe^{t})\frac{d}{dt}\big(e_{\lambda}(xe^{t})e_{\lambda}^{-1}(x)\big)&=(1+\lambda xe^{t})\sum_{n=1}^{\infty}\mathrm{Bel}_{n,\lambda}(x)\frac{t^{n-1}}{(n-1)!} \label{41} \\
    &=\sum_{n=0}^{\infty}\mathrm{Bel}_{n+1,\lambda}(x)\frac{t^{n}}{n!}+\lambda x\sum_{n=0}^{\infty}\bigg(\sum_{l=0}^{n}\binom{n}{l}\mathrm{Bel}_{l+1,\lambda}(x)\bigg)\frac{t^{n}}{n!}\nonumber.
\end{align}
Therefore, by comparing the coefficients on both sides of \eqref{40} and \eqref{41}, we obtain the following theorem.
\begin{theorem}
    For $n\ge 0$, we have
    \begin{displaymath}
        \mathrm{Bel}_{n+1,\lambda}(x)=x\sum_{l=0}^{n}\binom{n}{l}\big(\mathrm{Bel}_{l,\lambda}(x)-\lambda\mathrm{Bel}_{l+1,\lambda}(x)\big).
    \end{displaymath}
\end{theorem}
We observe that
\begin{align}
 \frac{d}{dx}\big(e_{\lambda}(xe^{t})\cdot e_{\lambda}^{-1}(x)\big)
=e^{t}e_{\lambda}(xe^t) e_{\lambda}^{-1}(x)(1+\lambda xe^{t}\big)^{-1}-e_{\lambda}(xe^{t})e_{\lambda}^{-1}(x)\cdot\frac{1}{1+\lambda x}. \label{42}
\end{align}
Thus, by \eqref{39} and \eqref{42}, we get
\begin{equation}
    \frac{d}{dt}\big(e_{\lambda}(xe^{t})\cdot e_{\lambda}^{-1}(x)\big)=x\bigg\{\frac{d}{dx}\big(e_{\lambda}(xe^{t})e_{\lambda}^{-1}(x)\big)+e_{\lambda}(xe^{t})e_{\lambda}^{-1}(x)\frac{1}{1+\lambda x}\bigg\}. \label{43}
\end{equation}
From \eqref{35} and \eqref{43}, we can derive the following identity
\begin{align}
    \sum_{n=0}^{\infty}\mathrm{Bel}_{n+1,\lambda}(x)\frac{t^{n}}{n!} &=\frac{d}{dt}\big(e_{\lambda}(xe^{t})\cdot e_{\lambda}^{-1}(x)\big) \nonumber\\
    &=x\bigg\{\frac{d}{dx}\big(e_{\lambda}(xe^{t})\cdot e_{\lambda}^{-1}(x)\big)+e_{\lambda}(xe^{t})\cdot e_{\lambda}^{-1}(x)\frac{1}{1+\lambda x}\bigg\}\nonumber \\
    &=x\bigg\{\sum_{n=0}^{\infty}\mathrm{Bel}_{n,\lambda}^{\prime}(x)+\frac{1}{1+\lambda x}\mathrm{Bel}_{n,\lambda}(x)\bigg\}\frac{t^{n}}{n!} \label{44}\\
    &=\frac{x}{1+\lambda x}\sum_{n=0}^{\infty}\bigg(\mathrm{Bel}_{n,\lambda}^{\prime}(x)+\mathrm{Bel}_{n,\lambda}(x)+\lambda x\mathrm{Bel}_{n,\lambda}^{\prime}(x)\bigg)\frac{t^{n}}{n!}, \nonumber
\end{align}
where $\displaystyle \mathrm{Bel}_{n,\lambda}^{\prime}(x)=\frac{d}{dx}\mathrm{Bel}_{n,\lambda}(x)\displaystyle$. \\
Comparing the coefficients on the both sides of \eqref{44}, we obtain the following theorem.
\begin{theorem}
    For $n\ge 0$, we have
    \begin{displaymath}
        \mathrm{Bel}_{n+1,\lambda}(x)=\frac{x}{1+\lambda x}\big(\mathrm{Bel}_{n,\lambda}^{\prime}(x)+\mathrm{Bel}_{n,\lambda}(x)+\lambda x\mathrm{Bel}_{n,\lambda}^{\prime}(x)\big),
    \end{displaymath}
where $\displaystyle \mathrm{Bel}_{n,\lambda}^{\prime}(x)=\frac{d}{dx}\mathrm{Bel}_{n,\lambda}(x)\displaystyle$.
\end{theorem}
By Theorem 9, we get
\begin{equation}
    \mathrm{Bel}_{n+1,\lambda}(x)=x\big(\mathrm{Bel}_{n,\lambda}^{\prime}(x)+\mathrm{Bel}_{n,\lambda}(x)+\lambda x \mathrm{Bel}_{n,\lambda}^{\prime}(x)-\lambda \mathrm{Bel}_{n+1,\lambda}(x)\big). \label{45}
\end{equation}
Thus, by Theorem 8 and \eqref{45}, we have
\begin{align}
    x\mathrm{Bel}_{n,\lambda}^{\prime}(x)&=\mathrm{Bel}_{n+1,\lambda}(x)-x\mathrm{Bel}_{n,\lambda}(x)-\lambda x^2\mathrm{Bel}_{n,\lambda}^{\prime}(x)+\lambda x\mathrm{Bel}_{n+1,\lambda}(x)\label{46} \\
    &=x\sum_{l=0}^{n}\binom{n}{l}\big(\mathrm{Bel}_{l,\lambda}(x)-\lambda\mathrm{Bel}_{l+1,\lambda}(x)\big)-x\mathrm{Bel}_{n,\lambda}(x)-\lambda x^2\mathrm{Bel}_{n,\lambda}^{\prime}(x)+\lambda x\mathrm{Bel}_{n+1,\lambda}(x)\nonumber \\
    &=x\sum_{l=0}^{n-1}\binom{n}{l}\big(\mathrm{Bel}_{l,\lambda}(x)-\lambda\mathrm{Bel}_{l+1,\lambda}(x)\big)-\lambda x^2\mathrm{Bel}_{n,\lambda}^{\prime}(x)\nonumber .
\end{align}
From \eqref{46}, we note that
\begin{equation}
    (1+\lambda x)\mathrm{Bel}_{n,\lambda}^{\prime}(x)=\sum_{l=0}^{n-1}\binom{n}{l}\big(\mathrm{Bel}_{l,\lambda}(x)-\lambda\mathrm{Bel}_{l+1,\lambda}(x)\big). \label{47}
\end{equation}
Therefore, by \eqref{47}, we obtain the following theorem.
\begin{theorem}
    For $n\in\mathbb{N}$, we have
    \begin{displaymath}
        \mathrm{Bel}_{n,\lambda}^{\prime}(x)=\frac{1}{1+\lambda x}\sum_{l=0}^{n-1}\binom{n}{l}\big(\mathrm{Bel}_{l,\lambda}(x)-\lambda\mathrm{Bel}_{l+1,\lambda}(x)\big).
    \end{displaymath}
\end{theorem}
From \eqref{8}, we note that
\begin{equation}
    \sum_{n=0}^{\infty} \mathrm{Bel}_{n,\lambda}(x)\frac{t^{n}}{n!}=e_{\lambda}(xe^{t})e_{\lambda}^{-1}(x)=\bigg(\frac{1+\lambda xe^{t}}{1+\lambda x}\bigg)^{\frac{1}{\lambda}}\label{48}
\end{equation}
\begin{displaymath}
    =\bigg(1+\frac{\lambda x}{1+\lambda x}(e^{t}-1)\bigg)^{\frac{1}{\lambda}}=\sum_{k=0}^{\infty}(1)_{k,\lambda}\bigg(\frac{x}{1+\lambda x}\bigg)^{k}\frac{1}{k!}(e^{t}-1)^{k}
\end{displaymath}
\begin{displaymath}
    =\sum_{k=0}^{\infty}(1)_{k,\lambda}\bigg(\frac{x}{1+\lambda x}\bigg)^{k}\sum_{n=k}^{\infty}S_{2}(n,k)\frac{t^{n}}{n!}=\sum_{n=0}^{\infty}\bigg(\sum_{k=0}^{n}(1)_{k,\lambda}\bigg(\frac{x}{1+\lambda x}\bigg)^{k}S_{2}(n,k)\bigg)\frac{t^{n}}{n!}.
\end{displaymath}
Therefore, by comparing the coefficients on both sides of \eqref{48}, we obtain the following theorem.
\begin{theorem}
    For $n\ge 0$, we have
    \begin{displaymath}
        \mathrm{Bel}_{n,\lambda}(x)=\sum_{k=0}^{n}(1)_{k,\lambda}\bigg(\frac{x}{1+\lambda x}\bigg)^{k}S_{2}(n,k).
    \end{displaymath}
    In particular, for $x=1$, we have
    \begin{displaymath}
        \mathrm{Bel}_{n,\lambda}=\sum_{k=0}^{n}(1)_{k,\lambda}\bigg(\frac{1}{1+\lambda}\bigg)^{k}S_{2}(n,k).
    \end{displaymath}
\end{theorem}
From \eqref{11}, we note that
\begin{equation}\begin{split}\label{49}
&\sum_{n=0}^{\infty}\int_{0}^{\infty}e^{-2x}\mathrm{bel}_{n,\lambda}(-x)dx\frac{t^{n}}{n!}=\int_{0}^{\infty}e^{-2x}e^{-x(e_{\lambda}(t)-1)}dx=\int_{0}^{\infty}e^{-x(e_{\lambda}(t)+1)}dx \\
&=\frac{1}{e_{\lambda}(t)+1}=\frac{1}{e_{\lambda}(t)-1}-\frac{2}{e_{\frac{\lambda}{2}}(2t)-1}=\frac{1}{t}\bigg(\frac{t}{e_{\lambda}(t)-1}-\frac{2t}{e_{\frac{\lambda}{2}}(2t)-1}\bigg)\nonumber\\
&=\frac{1}{t}\sum_{n=0}^{\infty}\big(\beta_{n,\lambda}-2^{n}\beta_{n,\frac{\lambda}{2}}\big)\frac{t^{n}}{n!}=\sum_{n=0}^{\infty}\bigg(\frac{\beta_{n+1,\lambda}-2^{n+1}\beta_{n+1,\frac{\lambda}{2}}}{n+1}\bigg)\frac{t^{n}}{n!}.\nonumber
\end{split}\end{equation}
Therefore, comparing the coefficients on both sides \eqref{49}, we obtain the following theorem.
\begin{theorem}
    For $n\ge 0$, we have
    \begin{displaymath}
        \int_{0}^{\infty}e^{-2x}\mathrm{bel}_{n,\lambda}(-x)dx=\frac{\beta_{n+1,\lambda}-2^{n+1}\beta_{n+1,\frac{\lambda}{2}}}{n+1}.
    \end{displaymath}
\end{theorem}
By \eqref{10} and \eqref{49}, we get
\begin{equation}\begin{split}\label{50}
&\sum_{n=0}^{\infty} \frac{\beta_{n+1,\lambda}-2^{n+1}\beta_{n+1,\frac{\lambda}{2}}}{n+1}\frac{t^{n}}{n!}=\int_{0}^{\infty}e^{-2x}e^{-x(e_{\lambda}(t)-1)}dt \\
&=\sum_{k=0}^{\infty}(-1)^{k}\int_{0}^{\infty}e^{-2x}x^{k}dx\frac {1}{k!}\big(e_{\lambda}(t)-1\big)^{k} =\sum_{k=0}^{\infty}(-1)^{k}\sum_{n=k}^{\infty}S_{2,\lambda}(n,k)\frac{t^{n}}{n!}\int_{0}^{\infty}e^{-2x}x^{k}dx \nonumber \\
 &=\sum_{n=0}^{\infty}\bigg(\sum_{k=0}^{n}(-1)^{k}S_{2,\lambda}(n,k)\frac{k!}{2^{k+1}}\bigg)\frac{t^{n}}{n!}.\nonumber
\end{split}\end{equation}
Therefore, by \eqref{50}, we obtain the following theorem.
\begin{theorem}
    For $n\ge 0$, we have
    \begin{displaymath}
        \sum_{k=0}^{n}(-1)^{k}2^{-k-1}k!S_{2,\lambda}(n,k)=\frac{1}{n+1}\big(\beta_{n+1,\lambda}-2^{n+1}\beta_{n+1,\frac{\lambda}{2}}\big).
    \end{displaymath}
\end{theorem}
For $p\in\mathbb{N}$, and $\mathrm{Re}(\delta)>0$, we have
\begin{align*}
    e_{\lambda}(x,\ \delta|-p)&=\sum_{n=0}^{\infty}\frac{(1)_{n,\lambda}x^{n}}{n!(n+\delta)^{-p}}=\sum_{n=0}^{\infty}\bigg[\sum_{k=0}^{p}\binom{p}{k}(1)_{n,\lambda}n^{k}\delta^{p-k}\bigg]\frac{x^{n}}{n!}\\
    &=\sum_{k=0}^{p}\binom{p}{k}\delta^{n-k}\bigg(\sum_{n=0}^{\infty}\frac{n^{k}}{n!}(1)_{n,\lambda}x^{n}\bigg)\\
    &=e_{\lambda}(x)\sum_{k=0}^{p}\binom{p}{k}\delta^{n-k}\bigg(e_{\lambda}^{-1}(x)\sum_{n=0}^{\infty}\frac{n^{k}}{n!}(1)_{n,\lambda}x^{n}\bigg) \\
    &=e_{\lambda}(x)\sum_{k=0}^{p}\binom{p}{k}\delta^{p-k}\mathrm{Bel}_{k,\lambda}.
\end{align*}
Therefore, with the help of Theorem 8 we obtain the following theorem.
\begin{theorem}
    For $\delta\in\mathbb{C}$ with $\mathrm{Re}(\delta)>0, \,\, any \ x\in\mathbb{C},\, and \, p\in\mathbb{N}$, we have
\begin{displaymath}
    e_{\lambda}(x,\ \delta|-p)=e_{\lambda}(x)\sum_{k=0}^{p}\binom{p}{k}\delta^{p-k}\mathrm{Bel}_{k,\lambda}(x).
\end{displaymath}
In particular, \,for $\delta=1$, we have
\begin{equation}\begin{split}
 e_{\lambda}(x,\ 1|-p)&=e_{\lambda}(x)\sum_{k=0}^{p}\binom{p}{k}\mathrm{Bel}_{k,\lambda}(x)\\
&=\frac{e_{\lambda}(x)}{x}\mathrm{Bel}_{p+1,\lambda}(x)+\lambda e_{\lambda}(x)\sum_{l=0}^{p}\binom{p}{l}\mathrm{Bel}_{l+1,\lambda}(x).\label{51}
\end{split}\end{equation}
\end{theorem}
Replacing $p$ by $p-1$ in \eqref{51}, we get
\begin{equation}
    e_{\lambda}(x,\ 1|1-p)=\frac{e_{\lambda}(x)}{x}\mathrm{Bel}_{p,\lambda}(x)+\lambda e_{\lambda}(x)\sum_{l=0}^{p-1}\binom{p-1}{l}\mathrm{Bel}_{l+1,\lambda}(x).\label{52}
\end{equation}
From \eqref{52}, we note that
\begin{align}
    e_{\lambda}(x,\ 1|1-p)xe_{\lambda}^{-1}(x)&=\mathrm{Bel}_{p,\lambda}(x)+\lambda x\sum_{l=0}^{p-1}\binom{p-1}{l}\mathrm{Bel}_{l+1,\lambda}(x) \label{53} \\
    &=(1+\lambda x)\mathrm{Bel}_{p,\lambda}(x)+\lambda x\sum_{l=0}^{p-2}\binom{p-1}{l}\mathrm{Bel}_{l+1,\lambda}(x).\nonumber
\end{align}
Thus, by \eqref{53}, we obtain
\begin{displaymath}
    \mathrm{Bel}_{p,\lambda}(x)=\frac{xe_{\lambda}^{-1}(x)}{1+\lambda x}e_{\lambda}(x,\ 1|1-p)-\frac{\lambda}{1+\lambda x}\sum_{l=0}^{p-2}\binom{p-1}{l}\mathrm{Bel}_{l+1,\lambda}(x), \,\,(p \geq 2).
\end{displaymath}
\section{Further Remarks}
For $\delta,s\in\mathbb{C}$ with $\mathrm{Re}(\delta)>0$ and $\mathrm{Re}(s)>1$, we observe that
\begin{equation}
 \int_{0}^{\infty}e_{\lambda}(tx, \delta|s)e^{-t}dt=\sum_{n=0}^{\infty}\frac{(1)_{n,\lambda}}{n!(n+\delta)^{s}}x^{n}\int_{0}^{\infty}t^{n}e^{-t}dt 
 =\sum_{n=0}^{\infty}\frac{(1)_{n,\lambda}}{(n+\delta)^{s}}x^{n} \label{54}.
\end{equation}
From \eqref{54}, we may consider the degenerate Lerch zeta function which is given by
\begin{equation}
    \Phi_{\lambda}(x,s,\delta)=\sum_{n=0}^{\infty}\frac{(1)_{n,\lambda}}{(n+\delta)^{s}}x^{n},\quad (\mathrm{Re}(\delta)>0).\label{55}
\end{equation}
Note that $\displaystyle\lim_{\lambda\rightarrow 0}\Phi_{\lambda}(x,s,\delta)=\Phi(x,s,\delta)\displaystyle$, where $\Phi(x,s,\delta)$ is Lerch zeta function defined by
\begin{displaymath}
    \Phi(x,s,\delta)=\sum_{n=0}^{\infty}\frac{x^{n}}{(n+\delta)^{s}},\quad(\mathrm{see}\ [1,18]).
\end{displaymath}
In particular, for $x=1$, we have
\begin{equation}
    \Phi_{\lambda}(1,s,\delta)=\sum_{n=0}^{\infty}\frac{(1)_{n,\lambda}}{(n+\delta)^{s}},\quad (\mathrm{Re}(\delta)>0).\label{56}
\end{equation}
In addition, Hurwitz zeta function is defined as
\begin{equation}
    \zeta(s,\delta)=\sum_{n=0}^{\infty}\frac{1}{(n+\delta)^{s}},\quad (\mathrm{Re}(\delta)>0). \label{57}
\end{equation}
In light of \eqref{57}, we denote $\Phi_{\lambda}(1,s,\delta)$ in \eqref{56} by $ \zeta_{\lambda}(s,\delta)$ and call it the degenerate Hurwitz zeta function as follows:\\
\begin{equation}
    \zeta_{\lambda}(s,\delta)=\sum_{n=0}^{\infty}\frac{(1)_{n,\lambda}}{(n+\delta)^{s}},\quad (\mathrm{Re}(\delta)>0). \label{58}
\end{equation}
Note that $\displaystyle\lim_{\lambda\rightarrow 0}\zeta_{\lambda}(s,\delta)=\zeta(s,\delta)\displaystyle$. \\
Let us take $\delta=1$. Then, by \eqref{58}, we get
\begin{equation}
    \zeta_{\lambda}(s,1)=\sum_{n=0}^{\infty}\frac{(1)_{n,\lambda}}{(n+1)^{s}}=\sum_{n=1}^{\infty}\frac{(1)_{n-1,\lambda}}{n^{s}} \label{59}.
\end{equation}
In view of \eqref{59}, we may consider the degenerate Riemann zeta function given by
\begin{equation}
    \zeta_{\lambda}(s) =\zeta_{\lambda}(s,1)=\sum_{n=1}^{\infty}\frac{(1)_{n-1,\lambda}}{n^{s}},\quad (\mathrm{Re}(s)>1). \label{60}
\end{equation}
By \eqref{3}, we get
\begin{align}
    \sum_{k=0}^{n-1}e_{\lambda}^{k}(t)&=\frac{1}{e_{\lambda}(t)-1}\big(e_{\lambda}^{n}(t)-1\big)=\frac{1}{t}\bigg\{\frac{t}{e_{\lambda}(t)-1}e_{\lambda}^{n}(t)-\frac{t}{e_{\lambda}(t)-1}\bigg\} \label{61} \\
    &=\frac{1}{t}\sum_{p=0}^{\infty}\big(\beta_{p,\lambda}(n)-\beta_{p,\lambda}\big)\frac{t^{p}}{p!}=\sum_{p=0}^{\infty}\bigg(\frac{\beta_{p+1,\lambda}(n)-\beta_{p+1,\lambda}}{p+1}\bigg)\frac{t^{p}}{p!}.\nonumber
\end{align}
On the other hand,
\begin{equation}
    \sum_{k=0}^{n-1}e_{\lambda}^{k}(t)=\sum_{k=0}^{n-1}\sum_{p=0}^{\infty}\frac{(k)_{p,\lambda}}{p!}t^{p}=\sum_{p=0}^{\infty}\bigg(\sum_{k=0}^{n-1}(k)_{p,\lambda}\bigg)\frac{t^{p}}{p!}. \label{62}
\end{equation}
From \eqref{61} and \eqref{62}, we have
\begin{equation}
    \sum_{k=0}^{n-1}(k)_{p,\lambda}=\frac{1}{p+1}\big\{\beta_{p+1,\lambda}(n)-\beta_{p+1,\lambda}\big\}, \label{63}
\end{equation}
where $n$ is a positive integer and $p$ is non-negative integer. By \eqref{63}, we easily get
\begin{displaymath}
    \sum_{k=0}^{n-1}(k)_{p,\lambda}=\frac{1}{p+1}\sum_{l=0}^{p}\binom{p+1}{l}(n)_{l,\lambda}\beta_{p+1-l,\lambda}.
\end{displaymath}


\begin{thebibliography}{9}
\bibitem{1}
 Aygunes, A. A.; Simsek, Y. \emph{Unification of multiple Lerch-zeta type functions,} Adv. Stud. Contemp. Math. (Kyungshang)  \textbf{21}  (2011),  no. 4, 367-373.
\bibitem{2}
Boyadzhiev, K. N.  \emph{Polyexponentials}, arXiv:0710.1332 [math.NA].
\bibitem{3}
Carlitz, L. \emph{Degenerate Stirling, Bernoulli and Eulerian numbers}, Utilitas Math. 15 (1979),  51-88.
\bibitem{4}
\emph{Comtet, L. Advanced combinatorics: the art of finite and infinite expansions (translated from the French by J.W. Nienhuys),} Dordrecht and Boston: Reidel, 1974.
\bibitem{5}
Dere, R.; Simsek, Y. \emph{Applications of umbral algebra to some special polynomials,} Adv. Stud. Contemp. Math. (Kyungshang)  \textbf{22}  (2012),  no. 3, 433-438.
\bibitem{6}
Dolgy, D.V.; Jang, G.-W.; Kim, T. \emph{A note on degenerate central factorial polynomials of the second kind}, Adv. Stud. Contemp. Math. (Kyungshang) 2019, \textbf{29} (1), 7-13.
\bibitem{7}
Dolgy, D. V.; Kim, T.; Kwon, H.-I.; Seo, J. J. \emph{Some identities for degenerate Euler numbers and polynomials arising from degenerate Bell polynomials,} Proc. Jangjeon Math. Soc.  \textbf{19}  (2016),  no. 3, 457-464.
\bibitem{8}
El-Desouky, B. S.; Mustafa, A. \emph{New results on higher-order Daehee and Bernoulli numbers and polynomials,} Adv. Difference Equ.  2016, Paper No. \textbf{32}, 21 pp.
\bibitem{9}
Hardy,G. H.  \emph{On the zeroes of certain classes of integral Taylor series. Part II.–-On the integral function formula and other similar functions,} Proc. London Math. Soc. (2)  \textbf{2}  (1905), 401-431.
\bibitem{10}
Hardy, G. H. \emph{On the zeroes certain classes of integral Taylor series. Part I.–-On the integral function formula,} Proc. London Math. Soc. (2)  \textbf{2}  (1905), 332-339.
\bibitem{11}
Kim, D. S.; Kim, T. \emph{A note on polyexponential and unipoly functions,} Russ. J.
Math. Phys.  \textbf{26} (2019), no. 1, 40-49.
\bibitem{12}
Jang, L.-C.; Kim, T.; Lee, D.-H.; Park, D.-W. \emph{An application of polylogarithms in the analogs of Genocchi numbers,} Notes Number Theory Discrete Math.  \textbf{7}  (2001),  no. 3, 65-69.
\bibitem{13}
Kim, T.; Kim, D. S.; Dolgy, D. V. \emph{On partially degenerate Bell numbers and polynomials,} Proc. Jangjeon Math. Soc.  \textbf{20}  (2017),  no. 3, 337-345.
\bibitem{14}
Kim, T.; Kim, D. S.; Kim, H.Y.; Kwon, J. \emph{Degenerate Stirling polynomials of the second kind and some applications,} Symmetry 2019, 11(8), 1046;
https://doi.org/10.3390/sym11081046.
\bibitem{15}
Kim, T.; Kim, D. S. \emph{Degenerate Laplace transform and degenerate gamma function,} Russ. J. Math. Phys. 2017, \textbf{24} (2), 241-248.
\bibitem{16}
Kim, T. \emph{A note on degenerate Stirling polynomials of the second kind,} Proc. Jangjeon Math. Soc.  \textbf{20}  (2017),  no. 3, 319-331.
\bibitem{17}
Kim, T.; Kim, D. S.; Lee,H.; Kwon, J. \emph{A note on some identities of new type
degenerate Bell polynomials,} Mathematics 2019, \textbf{7} (11), 1086
\bibitem{18}
Kim, T.; Rim, S.-H.; Simsek, Y.; Kim, D. \emph{On the analogs of Bernoulli and
Euler numbers, related identities and zeta and $L$ –functions,} J. Korean Math.
Soc.  \textbf{45} (2008),  no. 2, 435-453.
\bibitem{19}
Lalit Mohan, U. \emph{On the degenerate Laplace transform-I,} Bulletin of Pure \&
Applied Sciences- Mathematics and Statistics 37e (2018), no. 1, 1-8.
\bibitem{20}
Rim, S.-H.; Pak, H. K.; Kwon, J. K.; Kim, T. G. \emph{Some identities of Bell
polynomials associated with $p$-adic integral on $\mathbb{Z}_{p}$,} J. Comput. Anal. Appl. \textbf{20}
(2016), no. 3, 437-446.
\bibitem{21}
Whittaker, E. T.; Watson, G. N. A course of modern analysis. \emph{An introduction  to the general theory of infinite processes and of analytic functions: with an
account of the principal transcendental functions,} Fourth edition. Reprinted
Cambridge University Press, New York 1962 vii+608 pp.
\bibitem{22}
  Zhang, W.; Lin, X. \emph{Identities involving trigonometric functions and Bernoulli numbers,} Appl. Math. Comput. \textbf{334} (2018), 288-294.
\end{thebibliography}
\end{document}